\documentclass[12pt]{article}
\usepackage[centertags]{amsmath}
\usepackage{amsfonts}
\usepackage{amssymb}
\usepackage{amsthm}

\theoremstyle{plain}
\newtheorem{thm}{Theorem}[section]

\newtheorem{lem}[thm]{Lemma}

\newtheorem{remark}[thm]{Remark}
\theoremstyle{definition}
\newtheorem{defn}[thm]{Definition}
\newtheorem{exam}[thm]{Example}

\theoremstyle{remark}
\numberwithin{equation}{section}

\newcommand{\beast}{\begin{eqnarray*}}
\newcommand{\eeast}{\end{eqnarray*}}

\title{ An inequality for regular near polygons}

\author{Paul Terwilliger
\thanks
{Department of Mathematics, University of Wisconsin-Madison, USA}
\and Chih-wen Weng
\thanks {Department of Applied
Mathematics, National Chiao Tung University, Taiwan R.O.C.}}

\date{November 25, 2003}

\begin{document}
\maketitle

\bibliographystyle{plain}

\begin{abstract}
Let $\Gamma$ denote a near polygon distance-regular graph with
diameter $d\geq 3,$ valency $k$ and intersection numbers $a_1>0,$
$c_2>1.$  Let $\theta_1$ denote the second largest eigenvalue of
$\Gamma$. We show
\begin{equation*}
\theta_1\leq \frac{k-a_1-c_2}{c_2-1}.
\end{equation*}
We show the following (i)--(iii) are equivalent. (i) Equality is
attained above; (ii) $\Gamma$ is $Q$-polynomial with respect to
$\theta_1$; (iii) $\Gamma$ is a dual polar graph or a Hamming
graph.

\bigskip

\bigskip

\end{abstract}


\section{Introduction}\label{s1}
Let $\Gamma$ denote a near polygon distance-regular graph with
diameter $d\geq 3$ (see Section~\ref{s2} for formal definitions).
Suppose the intersection numbers $a_1>0$ and $c_2>1.$ It was shown
by Brouwer, Cohen and Neumaier that if $\Gamma$ has
classical parameters $(d, q, 0, \beta)$ then $\Gamma$ is a Hamming
graph or a dual polar graph \cite[Theorem~9.4.4]{bcn}. The same
conclusion was obtained by the second author under the assumption
that $\Gamma$ is $Q$-polynomial and has diameter $d\geq
4$ \cite[Corollary~5.7]{w:97}. Let
$\theta_0>\theta_1>\cdots>\theta_d$ denote the eigenvalues of
$\Gamma.$  It is known that $\theta_0=k,$ where $k$ denotes
the valency of $\Gamma.$  By \cite[Proposition~4.4.6(i)]{bcn},
$$\theta_d\geq -\frac{k}{a_1+1},$$
with equality if and only if  $\Gamma$ is a near $2d$-gon. We now
state our result.

\begin{thm}\label{main}
Let $\Gamma$ denote a near polygon distance-regular graph with
diameter $d\geq 3$,
valency $k$,
and intersection numbers $a_1>0,$ $c_2>1.$  Let $\theta_1$ denote
the second largest eigenvalue of $\Gamma$. Then
\begin{equation}\label{e4.1}
\theta_1\leq \frac{k-a_1-c_2}{c_2-1}.
\end{equation}
Moreover, the following (i)--(iii) are equivalent.
\begin{enumerate}
\item[(i)]  Equality  is attained in (\ref{e4.1}); \item[(ii)]
$\Gamma$ is $Q$-polynomial with respect to $\theta_1$;
\item[(iii)] $\Gamma$ is a dual polar graph or a Hamming graph.
\end{enumerate}
\bigskip

\end{thm}


\section{Preliminaries}\label{s2}
In this section we review some definitions and basic concepts. See
the books by Bannai and Ito \cite{BanIto} or Brouwer, Cohen, and
Neumaier \cite{bcn} for more  background information.
\bigskip

Let $\Gamma=(X,R)$ denote a finite, undirected, connected graph
without loops or multiple edges, with vertex set $X$, edge set
$R$, path-length distance function $\partial$ and diameter
$d:={\rm max}\{\partial(x,y)|x,y\in X\}.$ For $x\in X$ and for all
integers $i,$ set
$$\Gamma_i(x):=\{y|y\in X, \partial(x, y)=i\}.$$
Let $k$ denote a nonnegative integer. We say $\Gamma$ is {\it
regular} with {\it valency} $k$ whenever $|\Gamma_1(x)|=k$ for all
$x\in X.$ Pick an integer $i$ $(0\leq i\leq d).$ For $x\in X$ and
for $y\in \Gamma_i(x),$ set
\begin{eqnarray}
B(x, y)&:=&\Gamma_1(x)\cap\Gamma_{i+1}(y),\\
A(x, y)&:=&\Gamma_1(x)\cap\Gamma_i(y),\\
C(x, y)&:=&\Gamma_1(x)\cap\Gamma_{i-1}(y).
\end{eqnarray}
The graph $\Gamma$ is said to be {\it distance-regular} whenever
for all integers $i$\ $(0\leq i\leq d),$ and for all $x, y\in X$
with $\partial(x, y)=i,$ the numbers
\begin{equation}
c_i:=|C(x, y)|, \ \ \ \ \ \ a_i:=|A(x, y)|,\ \ \ \ \ \ b_i:=|B(x,
y)|
\end{equation}
are independent of $x$ and $y.$ We call the  $c_i,$ $a_i,$ $b_i$
the {\it intersection numbers} of $\Gamma.$ We
observe $c_0=0,$ $a_0=0,$ $b_d=0$ and $c_1=1.$  For the rest of
this paper we assume $\Gamma$ is distance-regular with diameter
$d\geq 3.$ We observe $\Gamma$ is regular with valence $k=b_0$ and
that
\begin{equation}\label{i-0}
c_i+a_i+b_i=k\ \ \ \ (0\leq i\leq d)
\end{equation}
\cite[p. 126]{bcn}.
\bigskip

We recall the Bose-Mesner algebra of $\Gamma$.
For $0 \leq i \leq d$ let
 $A_i$
 denote the matrix in $\hbox{Mat}_X(\mathbb{R})$ which
has $xy$ entry
\begin{equation}
\nonumber
(A_i)_{xy}=\begin{cases}
    1  & \text{if $\partial (x,y)=i$}\\
    0  & \text{if $\partial (x,y)\not=i$}
    \end{cases}\ \ \ \ \ (x,y\in X).
\end{equation}
We call $A_i$ the {\it $i$th distance matrix} of $\Gamma.$
Observe (ai) $A_0=I$; (aii) $\sum_{i=0}^d A_i=J$;  (aiii)
$A^t_i=A_i$ $(0\leq i\leq d)$, (aiv) $A_iA_j=\sum_{h=0}^d
p^h_{ij}A_h$ $(0\leq i,j \leq d)$, where $I$ denotes the identity
matrix and $J$ denotes the all ones matrix. We abbreviate $A:=A_1$
and call this the {\it adjacency matrix} of $\Gamma$. Let
$\mathbf{M}$ denote the subalgebra of ${\rm Mat}_X(\mathbb{R})$
generated by $A$. Using (ai)--(aiv)  we find $A_0, A_1, \cdots,
A_d$ form a basis of $\mathbf{M}$. We call $\mathbf{M}$ the {\it
Bose-Mesner algebra} of $\Gamma.$ By \cite[p. 59, p. 64]{BanIto},
$\mathbf{M}$ has a second basis $E_0$, $E_1$, $\cdots$, $E_d$ such
that (ei) $E_0=|X|^{-1}J$; (eii) $\sum_{i=0}^d E_i=I$; (eiii)
$E_i^t=E_i$ $(0\leq i\leq d)$; (eiv) $E_iE_j=\delta_{ij}E_i $
$(0\leq i,j\leq d)$. We call $E_0,$ $E_1,$ $\cdots$, $E_d$  the
{\it primitive idempotents} for $\Gamma$.
 Since $E_0,$ $E_1,$ $\cdots$, $E_d$
form a basis for
 $\mathbf{M}$
 there exist real scalars
$\theta_0,\theta_1,\cdots,\theta_d$ such that
  $A=\sum_{i=0}^d\theta_iE_i$.
By this and (eiv) we find $AE_i=\theta_i E_i$\ $(0\leq i\leq d)$.
Observe $\theta_0$, $\theta_1$, $\cdots$, $\theta_d$ are mutually
distinct since $A$ generates $\mathbf{M}.$ We assume the $E_i$ are
indexed so that $\theta_0>\theta_1>\cdots>\theta_d.$  We call
$\theta_i$ the {\it eigenvalue} of $\Gamma$ corresponding to
$E_i.$ By \cite[p. 197]{BanIto} we have $\theta_0=k$ and $-k\leq
\theta_i\leq k$\ $(0\leq i\leq d).$ We call $\theta_0$ the {\it
trivial} eigenvalue.

\bigskip

Let $\theta$ denote an eigenvalue of $\Gamma$ and let $E$ denote
the corresponding primitive idempotent. Since $E\in \mathbf{M},$
there exist
real numbers  $\sigma_0,$ $\sigma_1,$ $\cdots,$ $\sigma_d$
such that
\begin{equation}\label{7-14}
E=m|X|^{-1}\sum\limits_{i=0}^d\sigma_iA_i,
\end{equation}
where $m={\rm rank}\ E.$ We have $\sigma_0=1$ and
\begin{equation}\label{7-11}
c_i\sigma_{i-1}+a_i\sigma_i+b_i\sigma_{i+1}=\theta\sigma_i\ \ \ \
\ (0\leq i\leq d),
\end{equation}
where $\sigma_{-1},$ $\sigma_{d+1}$ denote indeterminates
\cite[p. 191]{BanIto}.
 The
sequence $\sigma_0,$ $\sigma_1,$ $\cdots,$ $\sigma_d$ is called
the {\it cosine sequence} associated with $\theta.$ Let
$\sigma_0,$ $\sigma_1,$ $\cdots,$ $\sigma_d$ denote the cosine
sequence associated with the eigenvalue $k.$ Comparing (\ref{i-0})
and (\ref{7-11}) we find $\sigma_i=1$\ \ $(0\leq i\leq d).$ By the
{\it trivial cosine sequence} of $\Gamma$ we mean the cosine
sequence associated with $k.$  Let $\theta$ denote
an eigenvalue of $\Gamma$ and let $\sigma_0, \sigma_1,\cdots,
\sigma_d$ denote the corresponding cosine sequence. By
(\ref{7-11}),
\begin{eqnarray}
\sigma_1&=&\theta k^{-1},\label{i2}\\
\sigma_2&=&\frac{\theta^2-a_1\theta-k}{kb_1}.\label{i3}
\end{eqnarray}
Combining (\ref{i2}) and (\ref{i3}) we find
\begin{equation}\label{0.1}
(\sigma_1-\sigma_2)b_1=(\theta+1)(\sigma_0-\sigma_1).
\end{equation}
\bigskip

Set $V=\mathbb{R}^X$ (column vectors). We define the inner product
$$\langle u, v\rangle=u^tv\ \ \ \ \ \ (u, v\in V).$$
For each $x\in X$ set
$$\hat x=(0, 0, \cdots, 1, 0, \cdots, 0)^t,$$
where
the $1$ is in coordinate $x.$ We observe $\lbrace{\hat x}|~x\in
X\rbrace$ is an orthonormal basis for $V.$ By (\ref{7-14}),
for $x,y \in X$ we have
\begin{equation}\label{i0}
\langle E\hat x, E\hat y\rangle=m|X|^{-1}\sigma_i,
\end{equation}
where $i=\partial(x, y).$

\bigskip

By a  {\it clique} in $\Gamma$ we mean a nonempty set consisting
of mutually adjacent vertices of $\Gamma.$ A given clique in
$\Gamma$ is said to be {\it maximal} whenever it is not properly
contained in a clique. The graph $\Gamma$ is said to be a {\it
near polygon} whenever
\begin{enumerate}
\item[(i)]  Each maximal clique has cardinality $a_1+2;$
\item[(ii)] For all maximal cliques $\ell$ and for all $x\in X,$
either
\begin{enumerate}
\item[(iia)] $\partial (x, y)=d$ for all $y\in\ell,$ or
\item[(iib)] there exists an integer $i$\ $(0\leq i\leq d-1)$ and
a unique $z\in \ell$ such that $\partial (x, z)=i$ and
$\partial(x, y)=i+1$ for all $y\in\ell-\{z\}.$
\end{enumerate}
\end{enumerate}
We give an alternate description of a near polygon. Let
$K_{1, 2,
1}$ denote the graph with $4$ vertices $s, x, y, s'$ such that
$\partial(s, x)=\partial(s, y)=\partial(x, y)=\partial(x,
s')=\partial(y, s')=1$ and $\partial(s, s')=2.$
 Then by \cite[Theorem 6.4.1]{bcn}
$\Gamma$ is a near polygon if and only if
both the following (i')-(ii') hold.
\begin{enumerate}
\item[(i')] $\Gamma$ does not contain an induced $K_{1, 2, 1}$
subgraph;
 \item[(ii')]  \begin{equation}\label{i-1} a_i=a_1c_i\ \ \
\ \ (0\leq i\leq d-1).
\end{equation}
\end{enumerate}
Assume $\Gamma$ is a near polygon. Then
\begin{equation}\label{8-11}
a_d\geq a_1c_d.
\end{equation}
Moreover $a_d=a_1c_d$ if and only if no maximal clique satisfies
(iia) above \cite[Theorem 6.4.1]{bcn}. In this case we call
$\Gamma$ a {\it near $2d$-gon}. Otherwise we call $\Gamma$ a {\it
near $(2d+1)$-gon}.
Assume $\Gamma$ is a near polygon.
The Hoffman bound states that
\begin{equation}\label{hoffman}
\theta_d\geq -\frac{k}{a_1+1},
\end{equation}
with equality if and only if $\Gamma$ is a near
$2d$-gon \cite[Proposition~4.4.6(i)]{bcn}.
\bigskip

\begin{defn}\label{db3.1}
Let $\Gamma$ denote a distance-regular graph with diameter $d\geq
3.$
We say $\Gamma$ has {\it classical parameters} $(d, q,
\alpha , \beta)$ whenever the intersection numbers are given by
\begin{equation}
c_i = {i\atopwithdelims [] 1}\biggl(1 + \alpha {i-1\atopwithdelims []
1}\biggr) \ \ \ \ \ \ \ (0\leq i \leq d), \label{b3.1}
\end{equation}
\begin{equation}
b_i = \biggl ({D\atopwithdelims [] 1} - {i\atopwithdelims [] 1}
  \biggr )\biggl (\beta - \alpha {i\atopwithdelims [] 1}\biggr )
  \ \ \ \ \ \ \ (0\leq i \leq d), \label{b3.2}
\end{equation}
where
\begin{equation}
{j\atopwithdelims [] 1} := 1 + q + q^2 + \cdots +q^{j-1}.
\label{b3.3}
\end{equation}
\end{defn}

We give two examples of near polygon distance-regular graphs with
classical parameters $(d,q,\alpha,\beta).$
\medskip

\begin{exam}\label{ex1}
{\bf The Hamming graph\  $H(d, n)$\ \  $(n\geq 2)$
           \cite{e:81}, \cite{n:85}, \cite{s:75},
           \cite{t:87}.}\medskip
\bigskip

\begin{enumerate}
\item [ ] $X$=all $d$-tuples of elements from the set $\{1,
2,\cdots ,n\},$

\item [ ] $xy\in R$ iff $x,y$ differ in exactly $1$ coordinate\ \
\
          $(x, y\in X),$
\item [ ] $q=1,$\ \ \ $\alpha=0,$\ \ \  $\beta=n-1,$

\item [ ] $c_i=i,$\ \ \ $b_i=(d-i)(n-1),$\ \ \ $a_i=(n-2)i$\ \ \ \
$(0\leq i\leq d),$

\item [ ] $\theta_i=\displaystyle (d-i)(n-1)-i\ \ \ \ (0\leq i\leq
d) .$
\end{enumerate}
\end{exam}

\begin{exam}\label{ex2} {\bf The Dual polar graphs \cite{c:82}, \cite{s:80}.}
\bigskip

\noindent Let $U$ denote a finite vector space with one of the
following non-degenerate forms:

$$
\begin{array}{ccccc}
\ \ \ \ {\rm name}\ \ \ \       &\ \ \ \ {\rm dim}(U)\ \ \ \    &\ \ \ \ {\rm field}\ \ \ \   &\ \ \ \ {\rm form}\ \ \  \        &\ \ \ \ \epsilon \ \ \ \ \\
\\
B_d(p^n)        &2d+1         &GF(p^n)       &{\rm quadratic}   &1        \\
\\
C_d(p^n)        &2d           &GF(p^n)       &{\rm symplectic}  &1        \\
\\
D_d(p^n)        &2d           &GF(p^n) &\frac{\displaystyle{{\rm
quadratic}}}{\displaystyle{({\rm Witt\ \ index}\ \ d)}}&0        \\
\\
~^2D_{d+1}(p^n) &2d+2         &GF(p^n) &\frac{\displaystyle{{\rm
quadratic}}}{\displaystyle{({\rm Witt\ \ index}\ \ d)}}&2        \\
\\
~^2A_{2d}(p^n)   &2d+1         &GF(p^{2n})     &{\rm Hermitean}   &\displaystyle{\frac{3}{2}}\\
\\
~^2A_{2d-1}(p^n) &2d           &GF(p^{2n})     &{\rm Hermitean}   &\displaystyle{\frac{1}{2}}\\
\end{array}
$$

A subspace of $U$ is called {\it isotropic} whenever the form
vanishes completely on that subspace. In each of the above cases,
the dimension of any maximal isotropic subspace is $d.$
\bigskip

\begin{enumerate}

\item [ ] $X=$ set all maximal isotropic subspaces of $U,$

\item [ ] $xy\in R$\ iff\ {\rm dim}$(x\cap y)$ $=$ $d-1$\ \ \ $(x,
y\in X),$

\item [ ] $\alpha=0,$\ \ \ $\beta=q^\epsilon,$

\item [ ] $c_i=\displaystyle{\frac {q^i-1}{q-1}},$\ \ \ \
$a_i=\displaystyle\frac{q^{i+\epsilon}-q^i-q^\epsilon+1}{q-1}$\ \
\ \ $(0\leq i\leq d),$

\item [ ]
$b_i=\displaystyle{\frac{q^{i+\epsilon}(q^{d-i}-1)}{q-1}}$\ \ \ \
$(0\leq i\leq d-1),$

\item [ ]
$\theta_i=\displaystyle{\frac{q^{d+\epsilon-i}-q^\epsilon-q^i+1}{q-1}
}$ \ \ \ \ $(0\leq i\leq d),$
\end{enumerate}

\noindent where
\begin{enumerate}

\item [ ] $q=p^n, p^n, p^n, p^n, p^{2n}, p^{2n}$  respectively.
\end{enumerate}
\end{exam}
\bigskip

The following three theorems will be used in the proof of our
results.
\bigskip

\begin{thm}\label{1.1}(\cite[Theorem 4.1]{t:95})
Let $\Gamma$ denote a distance-regular graph with diameter $d\geq
3,$ and let $q$ denote a real number at least $1.$ Then the
following conditions (i), (ii) are equivalent.
\begin{enumerate}
\item[(i)] $\Gamma$ has a nontrivial cosine sequence $\sigma_0,$
$\sigma_1,$ $\cdots ,$ $\sigma_d$ such that
$\sigma_{i-1}-q\sigma_i$ is independent of $i$\ \ $(1\leq i\leq
d).$ \item[(ii)] The intersection numbers of $\Gamma$ are such
that $qc_i-b_i-q(qc_{i-1}-b_{i-1})$ is independent of $i$\ \
$(1\leq i\leq d).$
\end{enumerate}
Furthermore, if (i), (ii) hold, then
\begin{equation}\label{8-11-2}
c_3\geq (c_2-q)(1+q+q^2).
\end{equation}

\end{thm}

\begin{thm}\label{1.2}(\cite[Theorem 4.2]{t:95})
Let $\Gamma$ denote a distance-regular graph with diameter $d\geq
3,$ and let $q$ denote a real number at least $1.$ Then the
following conditions (i), (ii) are equivalent.
\begin{enumerate}
\item[(i)] Statements (i), (ii) hold in Theorem~\ref{1.1}, and
$c_3=(c_2-q)(1+q+q^2).$ \item[(ii)] There exists $\alpha,\
\beta\in \mathbb{R}$ such that $\Gamma$ has classical parameters
$(d, q, \alpha, \beta).$
\end{enumerate}

\end{thm}

\begin{thm}\label{1.3} (\cite[Theorem 9.4.4]{bcn})
Let $\Gamma$ denote a distance-regular graph with diameter $d\geq 3$
with classical parameters $(d, q, 0, \beta).$ Assume the
intersection numbers $a_1>0$ and $c_2>1.$ Suppose $\Gamma$ is a
near polygon. Then $\Gamma$ is a dual polar graph or a
Hamming graph.
\end{thm}

\section{The inequality}

In this section we obtain the inequality in Theorem~\ref{main}.

\begin{lem}\label{2.1}
Let $\Gamma$ denote a near polygon distance-regular graph with
diameter $d\geq 3$,
valency $k$, and intersection numbers $a_1>0,$ $c_2>1.$ Let
$\theta_1$ denote the second largest eigenvalue of $\Gamma.$ Then
\begin{equation}\label{e0.01}
\theta_1\leq \frac{k-a_1-c_2}{c_2-1}.
\end{equation}
\end{lem}

\begin{proof} Abbreviate $E=E_1.$ Let
$\sigma_0, \sigma_1,\cdots, \sigma_d$ denote the cosine sequence
associated with $\theta_1.$ Fix any two vertices $x, y\in X$ with
$\partial(x, y)=2.$ We consider the vectors
\begin{eqnarray}
u&=&\sum\limits_{z\in A(x, y)} E\hat z -\sum\limits_{w\in A(y, x)}
E\hat w,   \label{8-6-1}\\
v&=&E\hat x-E\hat y.\label{8-6-2}
\end{eqnarray}
By the Cauchy-Schwartz inequality,
\begin{equation}\label{e1}
\| u \|^2
      \| v\|^2\geq    \langle u,\ v \rangle^2.
\end{equation}
We compute the terms in (\ref{e1}). Using (\ref{i0}),
(\ref{8-6-1}), (\ref{8-6-2}) we find
\begin{eqnarray}
 \| v \|^2&=&2m|X|^{-1}(\sigma_0-\sigma_2),\label{7-14-1}\\
 \langle u, v
\rangle&=&2m a_2|X|^{-1}(\sigma_1-\sigma_2).\label{7-14-2}
\end{eqnarray}
We now compute $\|u\|^2.$ To do this we first discuss the
distances between vertices in $A(x, y)$ and vertices in $A(y, x).$
We claim that for all $z\in A(x, y),$ $z$ is adjacent to $c_2-1$
vertices in $A(y, x)$ and is at distance $2$ from the remaining
$a_2-c_2+1$ vertices in $A(y, x).$ To see this fix $z\in A(x, y).$
Then $\ell:=A(x, z)\cup\{x, z\}$ is a maximal
clique; hence there exists a
unique vertex $s\in \ell$\ with $\partial(s, y)=1.$ That is $s\in
C(x, y)\cap C(z, y).$ Observe $|C(x, y)\cap C(z, y)|=1$, since any
other $s'\in C(x, y)\cap C(z, y)$ will cause either $xss'y$
or $sxzs'$ to be a $K_{1, 2, 1}$ subgraph. Hence there are $c_2-1$
vertices in $C(z, y)\cap A(y, x).$
 Observe for $w\in A(y, x)$ we have
$\partial(w, x)=2$ and $\partial(w, s)\leq 2$ so $\partial (w,
z)\leq 2.$ We have now proved the claim. Using the claim and
applying (\ref{i0}) we find
\begin{eqnarray}\label{7-14-3}
\| u\|^2\nonumber&=& \|\sum\limits_{z\in A(x, y)}E\hat z\|^2+\|
     \sum\limits_{w\in A(y, x)}E \hat w\|^2-2\langle
     \sum\limits_{z\in A(x, y)}E\hat z, \sum\limits_{w\in A(y, x)} E \hat
     w\rangle\nonumber\\
 &=&
2ma_2|X|^{-1}(\sigma_0+(a_1-c_2)\sigma_1+(c_2-a_1-1)\sigma_2).
\end{eqnarray}
Evaluating (\ref{e1}) using (\ref{7-14-1})--(\ref{7-14-3}) we
routinely find
\begin{equation}\label{e3}
(\sigma_0+(a_1-c_2)\sigma_1+(c_2-a_1-1)\sigma_2)(\sigma_0-\sigma_2)
\geq a_2(\sigma_1-\sigma_2)^2.
\end{equation}
Evaluating (\ref{e3}) using (\ref{i2}), (\ref{i3}), (\ref{i-1}) we
obtain
\begin{equation}\label{8-6-3}
(\theta_1-k)^2(\theta_1(a_1+1)+k)(k-\theta_1(c_2-1)-a_1-c_2)\geq
0.
\end{equation}
Clearly
$(\theta_1-k)^2>0$.
By (\ref{hoffman}) and since
$\theta_1>\theta_d$ we find $\theta_1(a_1+1)+k>0.$ Evaluating
(\ref{8-6-3}) using these comments we find
\begin{equation*}
k-\theta_1(c_2-1)-a_1-c_2\geq 0
\end{equation*}
and  (\ref{e0.01}) follows.
\end{proof}
\bigskip

\begin{remark}\label{8-6-5}
Referring to Example~\ref{ex1} and Example~\ref{ex2}, the eigenvalue
$\theta_1$ satisfies (\ref{e0.01}) with equality.
\end{remark}

We comment on the proof of Lemma~\ref{2.1}.

\begin{lem}\label{2.1b}
 With the notation of Lemma~\ref{2.1}, the following
(i)--(iii) are equivalent.
\begin{enumerate}
\item [(i)] Equality is attained in (\ref{e0.01}). \item[(ii)] For
all $x, y\in X$ such that $\partial(x, y)=2,$
\begin{equation}\label{e0.1} \sum\limits_{z\in A(x, y)}
E\hat z-\sum\limits_{w\in A(y, x)} E\hat w\in {\rm Span}(E\hat
x-E\hat y).
\end{equation}
\item[(iii)] There exist $x, y\in X$ such that $\partial(x, y)=2$
and
\begin{equation}\label{e0.1b} \sum\limits_{z\in A(x, y)}
E\hat z-\sum\limits_{w\in A(y, x)} E\hat w\in {\rm Span}(E\hat
x-E\hat y).
\end{equation}
Here $E=E_1.$
\end{enumerate}
\end{lem}
\begin{proof} Observe from the proof of Lemma~\ref{2.1} that
equality is attained
in (\ref{e0.01}) if and only if
equality is attained  in (\ref{e1}).
 We claim $v\not=0$. This will follow from
(\ref{7-14-1}) provided we can show $\sigma_0\not=\sigma_2$.
Suppose $\sigma_0=\sigma_2$.
Setting $\theta=\theta_1$ and $\sigma_2=\sigma_0$ in
(\ref{0.1}) and simplifying the result we find
$\theta_1=-b_1-1$. This is inconsistent with
(\ref{hoffman}) and  $\theta_1>\theta_d$.
We have now shown
 $\sigma_0\not=\sigma_2$ and it follows $v\not=0$.
We now see that equality is attained in (\ref{e1})
 if and only if $u\in {\rm Span}(v).$ The result follows.
\end{proof}

\section{The case of equality}

In this section we consider the case of equality in
(\ref{e0.01}).

\begin{lem}\label{8-14}
Let $\Gamma$ denote a near polygon distance-regular graph with
diameter $d\geq 3$ and intersection numbers
$a_1>0,$ $c_2>1.$  Let $\theta_1$ denote
the second largest eigenvalue of $\Gamma$ and let $\sigma_0,$
$\sigma_1,$ $\cdots,$ $\sigma_d$ denote the corresponding cosine
sequence. Suppose
equality holds  in
(\ref{e0.01}).
Then $\sigma_{i-1}-q\sigma_i$ is independent of $i$\ $(1\leq i\leq
d),$ where $q=c_2-1.$
\end{lem}

\begin{proof} Setting $c_2=q+1$
in (\ref{e0.01}) and using $k-a_1-1=b_1$ we find
  $\theta_1+1=b_1q^{-1}.$ In particular $\theta_1\not=-1.$
Observe
$\sigma_1\not=\sigma_2;$ otherwise $\sigma_0=\sigma_1$ by
(\ref{0.1}) forcing $\theta_1=k$ by (\ref{i2}), a contradiction.
Evaluating (\ref{0.1}) using  $\theta_1+1=b_1q^{-1}$ we find
\begin{equation}\label{8-1}
\frac{\sigma_0-\sigma_1}{\sigma_1-\sigma_2}=q.
\end{equation}
 Fix two vertices
$x, y\in X$ with $\partial(x, y)=2.$ Abbreviate $E=E_1.$ By
Lemma~\ref{2.1b} there exists $\lambda\in \mathbb{R}$ such that
\begin{equation}\label{e1.1}
\sum\limits_{z\in A(x, y)} E\hat z-\sum\limits_{w\in A(y, x)}
E\hat w=\lambda(E\hat x-E\hat y).
\end{equation}
Fix an integer $i$\ \ $(1\leq i\leq d-1)$ and pick
$u\in
X$ with $\partial(u, x)=i-1$ and $\partial(u, y)=i+1.$ Taking the
inner product of $E\hat u$ with both sides of (\ref{e1.1}),
\begin{equation}\label{e1.2}
a_2(\sigma_i-\sigma_{i+1})=\lambda (\sigma_{i-1}-\sigma_{i+1}).
\end{equation}
Setting $i=1$ in (\ref{e1.2}) we find
$a_2(\sigma_1-\sigma_2)=\lambda(\sigma_0-\sigma_2).$ From
(\ref{8-1}) we find $\sigma_0-\sigma_2=(\sigma_1-\sigma_2)(1+q).$
By these comments
$\lambda=a_2/(q+1).$
Evaluating (\ref{e1.2}) using this we find
$$\sigma_{i-1}-q\sigma_i=\sigma_i-q\sigma_{i+1}\ \ \ (1\leq i\leq
d-1).$$ From this we find $\sigma_{i-1}-q\sigma_i$ is independent
of $i$ for $1\leq i\leq d.$

\end{proof}

\begin{lem}\label{8-1-1}
Let $\Gamma$ denote a near polygon distance-regular graph with
$d\geq 3$ and intersection numbers $a_1>0,$ $c_2>1.$ Let
$\theta_1$ denote the second largest eigenvalue of $\Gamma$
and assume
 equality holds in
(\ref{e0.01}).
Then $\Gamma$ has
classical parameters $(d, q, 0, \beta).$
\end{lem}

\begin{proof} Let the scalar $q$ be as in Lemma~\ref{8-14}.
By Lemma~\ref{8-14}
we have Theorem~\ref{1.1}(i) and hence
Theorem~\ref{1.1}(ii). Applying Theorem~\ref{1.1}(ii) with $i=2,
3$ we find
 \begin{equation}\label{7-9}
qc_2-b_2-q(qc_1-b_1)=qc_3-b_3-q(qc_2-b_2).
 \end{equation}
 Simplifying (\ref{7-9}) using (\ref{i-0}) and $c_2=q+1,$ $a_2=a_1c_2$
  we obtain
\begin{equation}\label{7-9-2}
(a_1+1+q)(1+q+q^2-c_3)=a_3-a_1c_3.
\end{equation}
By (\ref{i-1}) we have $a_3=a_1c_3$ if $d>3,$ and by (\ref{8-11})
we have $a_3\geq a_1c_3$ if $d=3.$
 In any case
$a_3\geq a_1c_3$ so the right-hand side of (\ref{7-9-2}) is
nonnegative. Also $a_1+1+q>0$ since $q=c_2-1.$ Evaluating
(\ref{7-9-2}) using these comments we find
\begin{equation}
 c_3\leq 1+q+q^2. \end{equation}
 By (\ref{8-11-2}) and using $c_2=1+q$ we find
 $c_3\geq 1+q+q^2.$ Now apparently
$c_3=1+q+q^2.$ We can now check that the assumption
$c_3=(c_2-q)(1+q+q^2)$ in Theorem~\ref{1.2}(i) holds.
Applying
 Theorem~\ref{1.2} we find there exist
real numbers $\alpha, \beta $ such that
$\Gamma$ has classical parameters $(d, q,
\alpha, \beta)$.
By
(\ref{b3.1}) we find $c_2=(1+q)(1+\alpha)$.
By the construction $c_2=q+1$. Comparing these equations
we find $\alpha=0$.
\end{proof}

\noindent {\bf Proof of Theorem~\ref{main}.} The inequality
(\ref{e4.1}) is from (\ref{e0.01}).

\noindent (i)$\Longrightarrow$(iii).  By Lemma~{\ref{8-1-1},
$\Gamma$ has classical parameters $(d, q, 0, \beta).$ By this and
Theorem~\ref{1.3} we find $\Gamma$ is a dual polar graph or a
Hamming graph.

\noindent (iii)$\Longrightarrow$(ii) This is immediate from
\cite[Corollary 8.5.3]{bcn}.

\noindent (ii)$\Longrightarrow$(i) Lemma~\ref{2.1b}(ii) holds by
\cite[Theorem~3.3]{t:95}, so Lemma~\ref{2.1b}(i) holds and the
result follows.


\bigskip
\bigskip

\noindent Paul Terwilliger \hfil\break Department of Mathematics
\hfil\break University of Wisconsin-Madison \hfil\break Van Vleck
Hall \hfil\break 480 Lincoln Drive \hfil\break Madison, WI
53706--1388 \hfil\break USA \hfil\break
 email: {\tt terwilli@math.wisc.edu}
\hfil\break
\bigskip

\noindent Chih-wen Weng
 \hfil\break Department of Applied Mathematics
\hfil\break National Chiao Tung University \hfil\break 1001 Ta
Hsueh Road, \hfil\break Hsinchu  \hfil\break Taiwan 30050
\hfil\break email: {\tt weng@math.nctu.edu.tw} \hfil\break

\end{document}